\documentclass[11pt]{article}

\usepackage{enumerate}
\usepackage{graphicx}
\usepackage{pst-node,pstricks}
\usepackage{enumerate}
\usepackage{pst-node,pstricks}

 \font\tensym=msbm10
\font\sevensym=msbm7 \font\fivesym=msbm5

 \font\tengoth=eufb10
\font\sevengoth=eufb7 \font\fivegoth=eufb5

\newfam\symfam
\textfont\symfam=\tensym \scriptfont\symfam=\sevensym
\scriptscriptfont\symfam=\fivesym
\def\sym{\fam\symfam\tensym}

\newfam\gothfam
\textfont\gothfam=\tengoth \scriptfont\gothfam=\sevengoth
\scriptscriptfont\gothfam=\fivegoth

\def\bb{\sym}

\def\Box{\fbox{\hspace{2mm}}}
\def\Comm{{\rm Com}}

\def\ss{\smallskip}

\def\ra{\rightarrow}

\def\be{\beta}

\def\vt{\vartheta}
\def\S{\Sigma}

\def\Z{{\bb Z}}
\def\M{{\rm \bf M}}
\def\Rep{{\rm Rep}}
\def\Rec{{\rm Rec}}
\def\Rat{{\rm Rat}}
\def\TN{{\bf TN}}

\newtheorem{theoreme}{Th\'eor\`eme}
\newtheorem{proposition}[theoreme]{Proposition}

\newtheorem{note}{Note}

\newtheorem{example}[note]{Example}

\begin{document}
 \markboth{\sc G\'erard Duchamp and Jean-Gabriel Luque}
{Lazard's elimination is recognizable}
\title{  Lazard's elimination (in traces)
is finite-state recognizable  }
\author{
G\'erard {\sc Duchamp}$^\dag$
and Jean-Gabriel {\sc Luque}$^\ddag$\\ \\ \\
$^\ddag$Laboratoire d'Informatique de Paris Nord Universit\'e Paris
13\\ (LIPN UMR 7030 du CNRS),\\ Avenue
J.-B. Cl\'ement 93430 Villetaneuse, France\\
Gerard.Duchamp@lipn.univ-paris13.fr\\ \\
$^\ddag$ Laboratoire d'informatique de l'Institut Gaspard Monge\\
(UMR 8049 IGM-LabInfo ),\\ Universit\'e de Marne-la-Vall\'ee,\\ 5,
Boulevard
Descartes,\\ 775454 Champs-sur-Marne, France\\
Jean-Gabriel.Luque@univ-mlv.fr}
 \maketitle

\begin{abstract} We prove that the codes issued from the elimination of any
subalphabet in a trace monoid are finite-state recognizable. This
implies in particular that the transitive fatorizations of the
trace monoids are recognizable by (boolean) finite-state automata.
\end{abstract}
\section{Introduction}

Sch\"utzenberger (\cite{Lo} Chapter 5) introduced the notion of a
factorization
 of a monoid $M$
\begin{equation}
M=\prod^{\ra}_{i\in I} M_i
\end{equation}
where $(M_i)_{i\in I}$ is a subfamily of submonoids of the given
monoid $M$. When $M=A^*$ is a free monoid, at
 the both ends of the chain,
one has complete factorizations like Lyndon and
Hall
 factorizations \cite{Reut} and the bisections $|I|=2$ \cite{DT}.\\ A
nice way to produce factorizations is to start with a bisection
$M=M_1M_2$
 and refine the factors
using a uniform process. Doing this, we could obtain a complete
factorization for
 every trace monoid
\cite{DK}. Trace monoids are defined as follows. Consider an
alphabet $\Sigma=\{x_1,\cdots,x_n\}$ and a commutation relation
$\vt$ ({\it i.e.} a reflexive and symmetric relation) on $\Sigma$.
The trace monoid $\M(\S,\vt)$ is the quotient
\begin{equation}
\M(\S,\vt)=\S^*/_{\equiv_\vt}
\end{equation}
where $\equiv_\vt$ is the congruence generated by the relators
$ab\equiv ba$ where $(a,b)\in\vt$.\\ Later on,  we adressed the
question of bisecting a trace
 monoid so that the left factor be generated by a subalphabet
 (Lazard bisection) and the right
factor be a trace monoid \cite{DL}. Doing so, we obtained a
complete description of the factors and
 graph-theoretical
criteria for the factorization. We conjectured that the trace
codes so
 obtained could be
recognized by finite-state automata \cite{DL}.

\ss In this paper, we prove that the answer to the conjecture is
positive. This will be a consequence of
 the more general
result that if a trace monoid $M(\S,\vt)$ is bisected as
\begin{equation}
M(\S,\vt)=L.M(B,\vt_B)
\end{equation}
with $B\subset \S$ and $\vt_B=\vt \cap (B\times B))$, then the
minimal generating
 set $\be(L)$ of $L$
is recognizable by a finite-state, effectively constructible
automaton. Here, we prove this fact and give the construction of
the automaton.

The paper is organised as follows:

In section 2, we recall basic notions related to trace monoids and
recognizability. In section 3, we prove that the left factor of a
Lazard bisection is a recognizable set and we describe the
construction of a deterministic automaton recognizing it in
section 4. To end with, we explain in section 5 how to construct a
deterministic automaton which recognizes the generating set of the
left factor of such a bisection.
\section{Trace Monoids}
Trace monoids were introduced by Cartier and Foata with the
purpose of  studying some combinatorial problems linked with
rearrangements (see \cite{CF}). Next, this notion has been studied
by Mazurkiewicz and many schools of Computer Sciences in the context of concurrent program schemes (see
\cite{Maz,Maz2}).

 Let $x\in\S$ be a letter and denote
$\Comm(x)$ the set of letters which commute with $x$
\begin{equation}
\Comm(x)=\{z|(x,z)\in \vt\}.
\end{equation}
In particular, one has $x\in\Comm(x)$. Let $w\in\M(\S,\vt)$ be a
trace, we will denote \begin{equation} TA(w)=\{x\in \S|w=ux,
u\in\M(\S,\vt)\} \end{equation} the {\it terminal alphabet} of
$w$.

As it is shown in \cite{DK},  Lazard elimination occurs in the
context of traces. Let $B$ be a subalphabet of $\S$ and
$\vt_B=\vt\cap (B\times B)$. The trace monoid splits into two
submonoids
\begin{equation}
\M(\S,\vt)=L.\M(B,\vt_B)
\end{equation}
where $L$ is the submonoid consisting in the traces whose terminal
alphabet is a subset of $\S\setminus B$. Furthermore the
decomposition is unique, which suggests that  the following
equality occurs in $\Z\langle\S,\vt\rangle=\Z[\M(\S,\vt)]$, the
algebra of series corresponding to $\Z[\M(\S,\vt)]$ \cite{D34}.
Thus,
\begin{equation}
\underline{\M(\S,\vt)}=\underline{L}.\underline{\M(B,\vt_B)}
\end{equation}
where $\underline S$ denotes the characteristic series of a subset
$S\subset \M(\S,\vt)$ {\it i.e.}
\begin{equation}
\underline S=\sum_{w\in S}w\in
\Z\langle\langle\S,\vt\rangle\rangle.
\end{equation}
 Let $\phi$ be the natural
surjection $\S^*\rightarrow \M(\S,\vt)$, the {\it set of the
representative words }  of a trace $t$ is defined as
$\Rep(t)=\phi^{-1}(t)$. We can extend this definition to trace
langages $\Rep(L)=\phi^{-1}(L)=\bigcup_{t\in L}\phi^{-1}(t)$. A
trace langage is said {\it recognizable} if and only if its
representative set is, and we say that an automaton recognizes $L$
if and only if it {\it recognizes} $\Rep(L)$.
\begin{example}\rm
Let $a$ and $b$ be two commuting letters, then the set
$\Rep(\{ab\})$ is recognized by the automaton
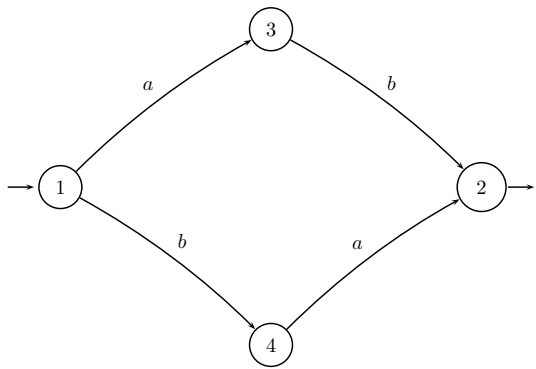
\begin{center}
\scalebox{0.7}{
\begin{pspicture}(8,-1)(10,8)
\cnode(5,3){15\pslinewidth}{b} \cnode(9,6){15\pslinewidth}{c}
\cnode(9,0){15\pslinewidth}{e} \cnode(13,3){17\pslinewidth}{d}
\ncarc{->}{b}{c}\Aput{$a$} \ncarc{->}{c}{d}\Aput{$b$}
\ncarc{->}{b}{e}\Aput{$b$} \ncarc{->}{e}{d}\Aput{$a$}
\psline{->}(4,3)(4.5,3) \psline{->}(13.5,3)(14,3)
 \rput(5,3){1}
\rput(13,3){2}\rput(9,6){3}\rput(9,0){$4$}
\end{pspicture}}
\end{center}
In fact, one can prove that a rational language is a set of
representatives (i.e. it is saturated w.r.t. the congruence
$\equiv_\theta$) if and only if the corresponding minimal
automaton shows complete squares as above.
\end{example}
We will denote $\Rec(\S,\vt)$ the set of  recognizable sets of
traces.

\section{Recognizing the left factor}
The $\Z$-rationality of the left factor $L$ is a direct
consequence of the unicity of the decomposition, which,  in term
of formal series, reads
\begin{equation}\label{}
\underline{\M(\S,\vt)}=\underline L.\underline{\M(B,\vt_B)}.
\end{equation}
where $\underline S$ denotes the $\Z$-characteristic series of the
set $S$ ({\it i.e.} $\underline S=\sum_{x\in S}x$). Indeed, by a
classical result due to Cartier and Foata (\cite{CF} Theorem 2.4)
the $\Z$-characteristic series of $\underline{\M(\S,\vt)}$ is
rational when the alphabet $\S$ is finite\footnote{The formula
holds also when the alphabet is infinite but the denominator is
then a series.} :

\begin{equation}
\underline{\M(\S,\vt)}={1\over\displaystyle\sum_{\{a_1,\dots,a_n\}\in{\bf
Cliques}(\S)}(-1)^na_1\cdots a_n}.
\end{equation}
where the sum at the denominator is taken over the set ${\bf
Cliques}(\S)$ of the cliques of $\S$ ({\it i.e.} commutative
sub-alphabets). Hence, one obtains the rational equality
\begin{equation}
\underline{L}= {1\over\displaystyle\sum_{\{a_1,\dots,a_n\}\in{\bf
Cliques}(\S)}(-1)^na_1\cdots a_n}\times \left(\displaystyle\sum_{\{b_1,\dots,b_n\}\in{\bf
Cliques}(B)}(-1)^nb_1\cdots
b_n\right)
\end{equation}
Nevertheless, this remark is not sufficient to show that $L$ is
 recognizable as a language. Furthermore, for traces, one
has the strict inclusion $\Rec(\S,\vt)\subset\Rat(\S,\vt)$.
To prove that $L$ is recognizable it suffices to find a
construction of $\Rep(L)$ using only recognizable operations.  For
each letter $x\in \S$, let $\TN_x$ be the set of representative
words of traces whose terminal alphabet does not contain $x$.
Remarking that $\Rep(L)$ is the representative set of the traces
whose terminal alphabet contains no letter of $B$, one has
\begin{equation}
\Rep(L)=\bigcap_{b\in B}\TN_b.
\end{equation}
Hence, $\Rep(L)$ is recognizable if each $\TN_b$ is. But, one can
easily verify that automaton ${\cal A}_{b}$:
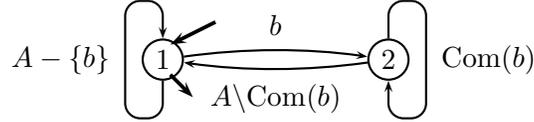
\begin{center}
\begin{pspicture}(5 ,1)
\cnode(1,.5){10\pslinewidth}{a} \cnode(4,.5){10\pslinewidth}{e}
\ncloop[angleA=270,angleB=90,loopsize=0.5,arm=0.5,linearc=.2]{->}{a}{a}\Aput{$A-\{b\}$}
\ncloop[angleA=90,angleB=270,loopsize=0.5,arm=0.5,linearc=.2]{->}{e}{e}\Aput{$\Comm(b)$}
\ncarc{->}{a}{e}\Aput{$b$}
\ncarc{->}{e}{a}\Aput{$A\backslash\Comm(b)$}
\pcline[linewidth=0.05]{<-}(1.1,.7)(1.7,1)
\pcline[linewidth=0.05]{->}(1.1,0.3)(1.4,-0.0) \rput(1,.5)1\rput(4
,.5)2
\end{pspicture}
\end{center}
recognizes $\TN_b$. Thus, we have the proposition
\begin{proposition}
 $L$ is a recognizable
submonoid of $M(\S,\vt).$
\end{proposition}

\section{A deterministic automaton for a terminal condition}
One can compute a deterministic automaton recognizing $L$
generalizing the construction of ${\cal A}_b$. We consider an
automaton ${\cal A}_B=(S_B,I_B,F_B, T_B)$ such that:
\begin{enumerate}
\item The set $S_B$ of its states  is the set of all the sub-alphabets
of $B$, \item There a unique initial state $I_B=\{\emptyset\}$,
\item There a unique final state $F_B=\{\emptyset\}=I_B$,
\item The transitions are $$T_B=\{(B',x,((B'\cup\{x\})\cap
\Comm(x)\cap B))\}_{B'\subset B, x\in \S}.$$
\end{enumerate}
One has
\begin{proposition}
The automaton ${\cal A}_B$ is a complete deterministic automaton
recognizing $\Rep(L)$.
\end{proposition}
{\bf Proof}  It is straightforward to see that such an automaton is
 complete and deterministic. Now, let us prove that it recognizes
$\Rep(L)$. As ${\cal A}_L$  is complete deterministic, for each
word $w=a_1\cdots a_n$ we can consider a state $s_w$ which is the
 state of ${\cal A}_B$ after  reading $w$. More
precisely, we can define $s_w$ as $s_w=s_n$ in the following chain
of transitions
\begin{equation}
(\emptyset,a_1,s_1), (s_1,a_2,s_2),\cdots, (s_{n-1},a_n,s_n=s_w).
\end{equation}
We first prove that if $w$ is a word then $s_w=TA(t_w)\cap B$
where $t_w$ denotes the trace admitting $w$ as representative
word. We use an induction process, considering as starting point:
$(\emptyset,x,\{x\}\cap B)$ where $x\in\S$. Let $w=a_1\cdots a_n$
be a word of length $n$, such that $s_w$ is the intersection
between $B$ and the terminal alphabet of trace $t_w$. Let
$a_{n+1}\in\sigma$ be an other letter. One has,
$(s_w,a_{n+1},s_{wa_{n+1}})\in T_B$. Hence, the set $s_{wa_{n+1}}$
is \begin{equation} s_{wa_{n+1}}=(s_w\cup \{a_{n+1}\})\cap
\Comm\{a_{n+1}\}\cap B=TA(t_wa_{n+1})\cap B=TA(t_{wa_{n+1}})\cap
B.
\end{equation}
This proves our assertion. Then, the set of words $w$ such that
$s_w=\emptyset$ is exactly the set of representative words of $L$.
$\Box$
\begin{example}\label{ex1}\rm
We consider the trace alphabet given by the following commutation
graph
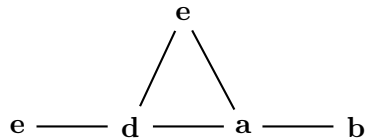
\begin{center}
\begin{pspicture}(10,2)
\psline(2.5,0)(7,0) \psline(4,0)(4.7,1.5) \psline(5.5,0)(4.7,1.5)
\rput(2.5,0){\pscirclebox[linestyle=none,fillstyle=solid]{\bf e}}
\rput(4,0){\pscirclebox[linestyle=none,fillstyle=solid]{\bf d}}
\rput(5.5,0){\pscirclebox[linestyle=none,fillstyle=solid]{\bf a}}
\rput(7,0){\pscirclebox[linestyle=none,fillstyle=solid]{\bf b}}
\rput(4.7,1.5){\pscirclebox[linestyle=none,fillstyle=solid]{\bf
e}}
\end{pspicture}
\end{center}
If we set $B=\{a,b\}$, then $L$ is recognized by the following
automaton (in the figure the only initial state and the only final
state is $\emptyset$):
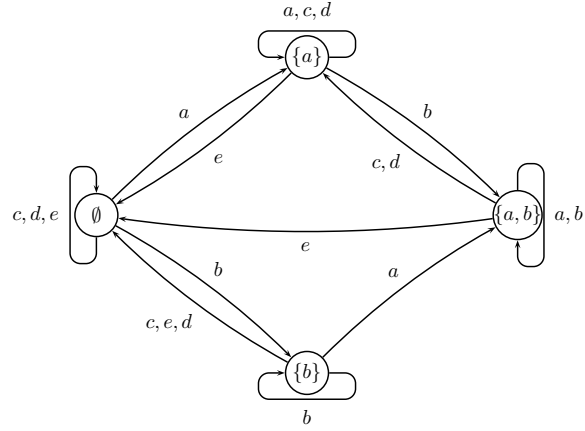
\begin{center}
\scalebox{0.7}{
\begin{pspicture}(8,-1)(10,8)
\cnode(5,3){15\pslinewidth}{b} \cnode(9,6){15\pslinewidth}{c}
\cnode(9,0){15\pslinewidth}{e} \cnode(13,3){17\pslinewidth}{d}
\ncarc{->}{b}{c}\Aput{$a$} \ncarc{->}{c}{d}\Aput{$b$}
\ncarc{->}{b}{e}\Aput{$b$} \ncarc{->}{e}{d}\Aput{$a$}
 \ncarc{->}{c}{b}\Aput{$e$} \ncarc{->}{e}{b}\Aput{$c,e,d$}
 \ncarc{->}{d}{c}\Aput{$c,d$}\ncarc{->}{d}{b}\Aput{$e$}
\ncloop[angleA=90,angleB=270,loopsize=0.5,arm=0.5,linearc=.2]{->}{d}{d}\Aput{$a,b$}
\ncloop[angleA=270,angleB=90,loopsize=0.5,arm=0.5,linearc=.2]{->}{b}{b}\Aput{$c,d,e$}
\ncloop[angleB=180,loopsize=0.5,arm=0.5,linearc=.2]{->}{c}{c}\Bput{$a,c,d$}
\ncloop[angleB=180,loopsize=-0.5,arm=0.5,linearc=.2]{->}{e}{e}\Aput{$b$}
\rput(5,3){$\emptyset$}
\rput(13,3){$\{a,b\}$}\rput(9,6){$\{a\}$}\rput(9,0){$\{b\}$}
\end{pspicture}}
\end{center}
\end{example}
\section{A deterministic automaton for  the generating set of the left factor}

Each  submonoid $M$ of a trace monoid has an unique {\it
generating set} which is the subset $G(M)=M\backslash
M^2$.\footnote{ The fact that $G(M)$ generates $M$ is
straightforward and the unicity comes from that the
$\Z$-characteristic series of $G(M)$ is the inverse of the
$\Z$-characteristic series of $M$ in $\Z\langle\langle
A\rangle\rangle$.}

In this section, we  prove that $G(L)$  is recognizable and we
construct an automaton $A_\beta$ which recognizes it. The
automaton $A_\beta$ is obtained from $A_B$ by adding two states
$F, H$, choosing $F$ as final state instead of ${\emptyset}$ and
modifying the transitions in such a way that if a letter of
$Z=A-B$ is read, the state reached belongs in ${F,H}$ and the
other states become unreachable.

More precisely, one considers the automaton ${\cal
A}_\beta=(S_\beta,I_\beta,F_\beta, T_\beta)$ obtained from the
automaton $A_B=(S_B,I_B,F_B,T_B)$ computed in the previous section
as follows:
\begin{enumerate}
\item The set of its states $S_\beta$, is the set of the sub-alphabets
of $B$ plus two states $F$ and $H$, \item There is a unique
initial state $I_\beta=\{\emptyset\}$,
\item There is a unique final state $F_\beta=\{F\}$,
\item The transitions are $$T_\beta=T_{B{\rightarrow} B}\cup
T_{B\rightarrow F} \cup T_{B\rightarrow H}\cup T_{F\rightarrow H}
\cup T_{H\rightarrow H}$$ where
\begin{enumerate}
\item $T_{B\rightarrow B}=\{(B',b,B'')\}_{B',B''\subset B, b\in B,(B',b,B'')\in T_B},$
\item $T_{B\rightarrow F}=\{(B',z,F)\}_{(B',z,\emptyset)\in T_B, B'\subset
B,B'\neq\emptyset, z\in Z},$
\item $T_{B\rightarrow H}=\{(B',z,H)\}_{(B',z,B'')\in
T_B, B',B''\not\in\{\emptyset, F,H\}, z\in Z},$
\item $T_{F\rightarrow H}=\{(F,x,H)\}_{x\in\S},\,$
\item and  $T_{H\rightarrow H}=\{(H,x,H)\}_{x\in\S}$.
\end{enumerate}
\end{enumerate}
\begin{proposition}
The automaton ${\cal A}_\beta$ recognizes $\Rep(G(L))$.
\end{proposition}
{\bf Proof}
 The automaton is almost the same as ${\cal A}_B$. As for
${\cal A}_B$, if a word of $B^*$ is read, the automaton is in the
state corresponding to its terminal alphabet. The difference
appears when a letter of $Z$ is read, if it is read from the
$\emptyset$ state the automaton goes to the state $F$. Consider
now a word $w=w'z$ with $w'\in B^+$, $z\in Z$. We denote
$\delta_w$ the state of the automaton after reading $w$ (this
definition makes sense as, like ${\cal A}_B$, ${\cal A}_\beta$ is
deterministic). Now, if $\{z\}=TA(w'z)$, then $(\delta_w',z,F)\in
T_{B\rightarrow F}$ which means that $w$ is recognized by ${\cal
A}_\beta$, otherwise
 $(\delta_w',z,H)\in T_{B\rightarrow H}$ and $w$ is not recognized by ${\cal
 A}_\beta$. Furthermore, for each $z\in Z$ and $b\in B$,
 $\delta_{w'zaw''}=H$ (for each $w',w''\in \S^*$). This ends the
 proof.$\Box$
\begin{example}\rm
Consider again the example (\ref{ex1}). Then, $\beta$ is
recognized by the automaton
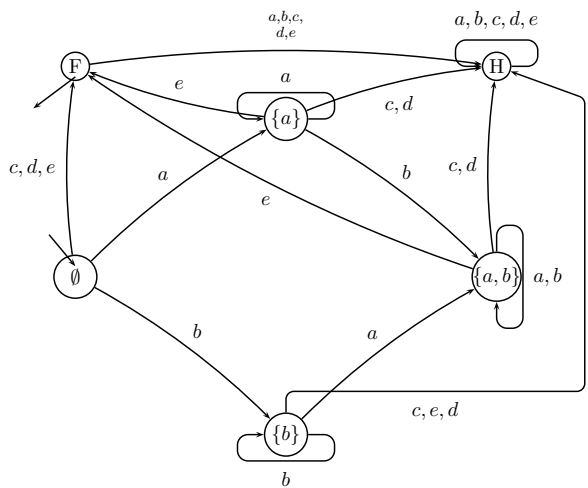
\begin{center}
\scalebox{0.7}{
\begin{pspicture}(8,-1)(10,10)
\cnode(5,3){15\pslinewidth}{b} \cnode(9,6){15\pslinewidth}{c}
\cnode(9,0){15\pslinewidth}{e} \cnode(13,3){17\pslinewidth}{d}
\cnode(13,7){10\pslinewidth}{H}\cnode(5,7){10\pslinewidth}{F}
 \ncarc{->}{b}{c}\Aput{$a$} \ncarc{->}{c}{d}\Aput{$b$}
\ncarc{->}{b}{e}\Aput{$b$} \ncarc{->}{e}{d}\Aput{$a$}
 \ncarc{->}{c}{H}\Bput{$c,d$} \ncarc{->}{c}{F}\Bput{$e$}
 \ncangles[angleA=90,armA=0.4cm,armB=1.5cm,angleB=-20,linearc=.15]{->}{e}{H}\Bput{$c,e,d$}
 \ncarc{->}{d}{H}\Aput{$c,d$}\ncarc{->}{d}{F}\Aput{$e$}
 \ncarc{->}{b}{F}\Aput{$c,d,e$}\ncarc{->}{F}{H}\Aput{$a,b,c,\atop d,e$}
\ncloop[angleA=90,angleB=270,loopsize=0.5,arm=0.5,linearc=.2]{->}{d}{d}\Aput{$a,b$}
\ncloop[angleB=180,loopsize=0.5,arm=0.5,linearc=.2]{->}{c}{c}\Bput{$a$}
\ncloop[angleB=180,loopsize=-0.5,arm=0.5,linearc=.2]{->}{e}{e}\Aput{$b$}
\ncloop[angleB=180,loopsize=0.5,arm=0.5,linearc=.2]{->}{H}{H}\Bput{$
{a,b,c, d,e}$} \rput(5,3){$\emptyset$}
\rput(13,3){$\{a,b\}$}\rput(9,6){$\{a\}$}\rput(9,0){$\{b\}$}
\rput(5,7)F\rput(13,7)H
\pcline{->}(5,6.8)(4.2,6.2)\pcline{<-}(5,3.2)(4.5,3.8)
\end{pspicture}}
\end{center}
\end{example}
\section{Conclusion}
The factorisations of free monoids (or in a more general setting
of a monoid constructed by generators and relations) is a relevant
topic in the context of the theory of codes \cite{BP}. Lazard
bisections, or more generally rational bisections \cite{DT}, play
a role in the construction of bases of free Lie algebras
\cite{Reut} and the study of circular codes \cite{BP,Reut}. A
natural question asks if it is possible to generalize these
properties to other monoids in particular when the free module
over these monoids can be endowed with a shuffle coproduct
\cite{DL2}. The results contained in the paper consist in a step
in the study of these problems
 for the trace monoids. The role played by
  the Lazard bisections in this context is not
 still completely known (see \cite{DK,DL} for some results).

\end{document}